\documentclass[12pt]{amsart}
\usepackage{epic}
\usepackage{eepic}
\usepackage{amsthm}
\usepackage{amscd}

\newcommand     {\comment}[1]   {}
\newcommand{\mute}[2] {}
\newcommand     {\printname}[1] {}
\newcommand{\labell}[1] {\label{#1}\printname{#1}}

\def \sss {\scriptstyle}

\DeclareMathOperator {\tDH} {DH}
\def    \M {{\mathcal{M}}}

\def    \J {{\mathcal{J}}}
\def \calM {{\mathcal M}}

\def \half {{\frac12}}

\def \t {{\mathfrak t}}
\def \g {{\mathfrak g}}

\def \Z {{\mathbb Z}}
\def \R {{\mathbb R}}
\def \C {{\mathbb C}}
\def \CP {{\mathbb C}{\mathbb P}}
\def \tCP {{\mathbb C}{\mathbb P}{^2}}
\def \bfZ {{\mathbf{Z}}}
\def \eps {\epsilon}
\def \inv {^{-1}}
\def \tomega {\widetilde{\omega}}
\def \tM {\widetilde{M}}
\def \tPhi {\widetilde{\Phi}}

\DeclareMathOperator {\GL} {GL}
\DeclareMathOperator {\PSL} {PSL}

\def \FS {{\text{FS}}}

\numberwithin{equation}{section} 
\newtheorem{Theorem}{Theorem}[section] 
\newtheorem{Lemma}{Lemma}[section] 
\theoremstyle{remark} 
\newtheorem{Remark}{Remark}[section] 
\newtheorem*{Remarks}{Remarks} 
\newtheorem*{Example}{Example}

\begin{document}

\title[Torus actions on equal symplectic blow-ups 
of $\mathbf{\text{CP}^2}$]{Circle and torus actions on equal symplectic blow-ups 
of $\mathbf{\text{CP}^2}$}

\author{Yael Karshon}
\address{Department of Mathematics, University of Toronto, 
Toronto, ON M5S 3G3, Canada.}
\email{karshon@math.toronto.edu}

\author{Liat Kessler}
\address{Courant Institute of Mathematical Sciences,
New York University, New York, NY 10012, U.S.A.}
\email{kessler@cims.nyu.edu}

\begin{abstract}
A manifold obtained by $k$ simultaneous symplectic blow-ups of $\CP^2$
of equal sizes $\epsilon$ 
(where the size of $\CP^1\subset\CP^2$ is one)
admits an effective two dimensional torus action if $k \leq 3$ and admits 
an effective circle action if $(k-1)\epsilon < 1$.
We show that these bounds are sharp if $1/\epsilon$ is an integer.
\end{abstract}

\maketitle

\section{Toric actions and circle actions in dimension four}
\labell{sec:actions}

\subsection*{Hamiltonian torus actions} 

Let a torus $T \cong (S^1)^k$ act on a compact connected
symplectic manifold $(M,\omega)$ of dimension $2n$
by symplectic transformations.
The action is \emph{Hamiltonian} if there exists
a moment map, that is, a map $$\Phi \colon M \to \t^* \cong \R^k$$
such that 
$$ d \Phi_j = - \iota(\xi_j) \omega $$ 
for all $j=1,\ldots,k$, where $\xi_1,\ldots,\xi_k$ are the vector fields 
that generate the torus action. If $H^1(M)=0$ then every
symplectic torus action is Hamiltonian. 
By the convexity theorem \cite{atiyah,GS:convexity},
the image of the moment map, 
$$\Delta:= \Phi(M),$$ 
is a convex polytope.
By the equivariant Darboux-Weinstein theorem \cite{W},
every $T$-fixed point has a neighborhood $U$
which is equivariantly symplectomorphic to a neighborhood
of the origin in $\C^n$ with $T$ acting linearly.
The components $\Phi^\xi = \left< \Phi , \xi \right>$, $\xi \in \g$,
of the moment map are perfect Morse-Bott functions.
See \cite{GS:convexity}.

From now on we assume that the action of $T$ is effective.

\subsection*{The Delzant theorem}
If $\dim T = \half \dim M$, the triple $(M,\omega,\Phi)$ 
is a \emph{symplectic toric manifold}, and the $T$-action is called 
\emph{toric}. 
By the Delzant theorem (\cite{delzant}; also see \cite{LT}), 
$(M,\omega,\Phi)$ is determined by $\Delta$ up to an equivariant 
symplectomorphism 
preserving $\Phi$. The inverse image under $\Phi$ of a vertex of $\Delta$ 
is a fixed point for the $T$-action, and the image of a $T$-fixed point 
is a vertex of $\Delta$.

A necessary condition for $\Delta$ to occur as the moment map image 
of a symplectic toric manifold is that it be a \emph{Delzant polytope}, 
meaning that the edges emanating from each vertex are generated 
by vectors $v_1,\ldots,v_n$ that span the lattice $\Z^n$.

\subsection*{Delzant's construction}
Given a Delzant polytope $\Delta$, Delzant constructs 
a symplectic toric manifold $(M_\Delta,\omega_\Delta,\Phi_{\Delta})$
whose moment map image is $\Delta$.  This manifold is a symplectic 
quotient of $\C^N$, where $N$ is the number of facets of $\Delta$,
with respect to a subgroup $K$ of $(S^1)^N$;
the toric $T$-action is through an isomorphism of $T$
with the quotient $(S^1)^N/K$.
The polytope $\Delta$ is then realized as the intersection
of the positive orthant $\R_+^N$ with an affine plane.
See \cite{delzant}, \cite{audin}, or \cite{g:moment}.

\subsection*{Toric actions on $\mathbf{CP^2}$}

An important special case of the Delzant construction
is the construction of $\CP^2$ as the quotient of the sphere
$S^5 \subset \C^3$ by the diagonal $S^1$-action, resulting in the
\emph{Fubini-Study} symplectic form $\omega_\FS$ on $\CP^2$.
Whenever we refer to $\CP^2$ as a symplectic manifold,
we assume that the 
symplectic form is $\omega_\FS$, normalized so that
$$\frac{1}{2\pi} \int_{\CP^1} \omega_\FS = 1.$$

The standard toric action on $\CP^2$ is
$(a,b) \cdot [z_0 : z_1 : z_2] = [z_0 : a z_1 : b z_2]$.
The moment map image is the triangle 
\begin{equation} \labell{Delzant triangle}
 \{ (x_1 , x_2) \in \R^2 \ | \ x_1 \geq 0 , x_2 \geq 0 , 
     x_1 + x_2 \leq 1 \}.
\end{equation}

By the Delzant theorem, every toric action on $\CP^2$ is obtained
from a symplectomorphism of $\CP^2$ with a symplectic toric 
manifold $M_\Delta$.  The second Betti number $\dim H_2(M_\Delta)$
is equal to the number of edges of $\Delta$ minus two;
this follows from Morse theory for the moment map.
So $\Delta$ must be a triangle.
It is easy to check that every Delzant triangle can be obtained from 
a standard one, \eqref{Delzant triangle}, by a transformation
$x \mapsto Ax+b$ where $A \in \GL(2,\Z)$ and $b \in \R^2$.
It follows that 

\begin{Lemma} \labell{toric is standard}
Every toric $T$-action on $\CP^2$ is equivariantly symplectomorphic 
to the standard action through an isomorphism of the torus $T$ 
with $(S^1)^2$.  
\end{Lemma}

\subsection*{Hamiltonian circle actions on compact symplectic four-manifolds}

Let the circle group $T=S^1$ act on a compact connected symplectic manifold 
$(M,\omega)$ of dimension 4 with moment map
$\Phi \colon M \to \t^* \cong \R$.
A connected component of the fixed point set is either an isolated 
fixed point $p \in M$ or a closed symplectic surface $F \subset M$
on which $\Phi$ is constant.  
Every interior fixed point (that is, a fixed point at which  $\Phi$
   is neither maximal nor minimal) is isolated.
Let $\bfZ_k \subset S^1$ 
be the cyclic subgroup of order $k$.  A connected component of the 
fixed point set of $\bfZ_k$
which is not fixed by any larger subgroup is a symplectic 2-sphere
$C \subset M$ on which $S^1$ acts by rotations of ``speed k";
we call it a \emph{$\bfZ_k$-sphere}.
A $\bfZ_k$-sphere contains two fixed points, $p$ and $q$, 
which are isolated fixed points in $M$.  
We note that $\frac{1}{k} \left| \Phi(p) - \Phi(q) \right| = 
\frac{1}{2\pi} \int_C \omega$.

To $(M,\omega,\Phi)$ we associate the following labeled graph.
To an isolated fixed point $p$ we associate a vertex $\langle p \rangle$,
labeled by the real number $\Phi(p)$.
To a $\bfZ_k$-sphere containing two fixed points $p$ and $q$
we associate an edge connecting the vertices $\langle p \rangle$
and $\langle q \rangle$ and labeled by the integer $k$.
To a two dimensional component $F$ of the fixed point set we associate
a vertex $\langle F \rangle$ labeled by two real numbers and one
integer: the moment map value $\Phi(F)$, the symplectic size
$\frac{1}{2\pi} \int_F \omega$, and the genus of the surface $F$.
The graph determines $(M,\omega,\Phi)$ up to an equivariant symplectomorphism
that preserves moment maps \cite{karshon}.

\subsection*{Circle actions on $\mathbf{CP^2}$}
On a compact symplectic four-manifold $(M,\omega)$,
if $\dim H^2(M;\R) \leq 3$ and $\dim H^1(M;\R) =0$,
then every Hamiltonian circle action on $(M,\omega)$
extends to a toric action.  This is proved in
\cite[Theorem 1]{max-torus} using techniques of \cite{karshon}.
In particular, 

\begin{Lemma} \labell{S1 is standard}
Every Hamiltonian $S^1$-action on $\CP^2$ is the composition 
of a homomorphism $S^1 \to (S^1)^2$ with the standard toric action 
on $\CP^2$.
\end{Lemma}

\section{Symplectic blow-ups}
\labell{sec:blowup}

\subsection*{Blow-up of $\mathbf{C^n}$}

Consider $\C^n \cong \R^{2n}$ with its standard symplectic form
$\sum d x_i \wedge d y_i$.  The standard symplectic blow-up of $\C^n$
of \emph{size} $r^2/2$ is obtained by removing
the open ball $B^{2n}(r)$ of radius $r$ about the origin 
and collapsing its boundary along the Hopf fibration 
$S^{2n-1} \to \CP^{n-1}$. 
The resulting space is the disjoint union of $\CP^{n-1}$
with the open set $\{ \| z \| > r \}$.
This space is naturally a smooth symplectic manifold such that
the symplectic form on $\{ \| z \| > r \}$ is induced from $\C^n$
and the symplectic form on $\CP^{n-1}$, pulled back to $S^{2n-1}$,
is also induced from $\C^n$.  
This can be shown through symplectic cutting \cite{lerman} 
or by an explicit formula on the complex blow-up \cite{GS:birational}.
For more details see \cite[Section 7.1]{MS:intro}.
The submanifold $\CP^{n-1}$ is called the \emph{exceptional divisor}.

\subsection*{Blow-ups of a symplectic manifold}

A blow-up of a $2n$-dimensional symplectic manifold $(M,\omega)$
is a new symplectic manifold that is constructed in the following way.
Let $\Omega \subset \C^n$ be an open subset that contains a ball 
about the origin of radius greater than $r$, and let 
$i \colon \Omega \to M$ be a symplectomorphism onto an open subset of $M$.
The standard blow-up in $\Omega$ of size $r^2/2$ transports to $M$
through $i$.

Let $(\tilde{M},\tilde{\omega})$ denote the resulting manifold.
Then 
\begin{equation} \labell{property1}
\dim H_2(\tilde{M}) = \dim H_2(M) + 1,
\end{equation}
and 
\begin{equation} \labell{property2}
\frac{1}{(2\pi)^n} \int_{\tilde{M}} \tomega^n
 \ = \  \frac{1}{(2\pi)^n} \int_M \omega^n \ - \  \epsilon^n
\end{equation}
where $\epsilon = r^2/2$ is the size of the blow-up. 

Similarly, $k$ \emph{simultaneous blow-ups} are obtained from embeddings
$i_1 \colon \Omega_1 \to M$, $\ldots,$ $i_k \colon \Omega_k \to M$
whose images are disjoint.

\subsection*{Notions of size}

The \emph{size}, or \emph{normalized symplectic area}, 
of a two dimensional symplectic manifold $(C,\omega)$
is $\frac{1}{2\pi} \int_C \omega$.
We also define the size of a ball of radius $r$ in $\C^n$
to be $r^2/2$, and the size of $(\CP^n,\omega)$ to be 
$\frac{1}{2\pi} \int_{\CP^1} \omega$.
These notions are compatible with the notion of the size
of a symplectic blow-up: in a symplectic blow-up of size $\epsilon$, 
the size of the resulting exceptional divisor is $\epsilon$.

\subsection*{Constraints on blow-ups}

The volume of a symplectic manifold gives constraints on the sizes of 
its symplectic blow-ups: the volume of an embedded ball cannot exceed 
the volume of the manifold.   
Sharper constraints are proved using holomorphic techniques.
For instance, by the volume constraint, any symplectic blow-up of $\CP^2$
must have size $< 1$.  By holomorphic constraints, any two simultaneous 
symplectic blow-ups of $\CP^2$ must have sizes whose sum is $<1$.  
See \cite{gromov}.  For further constraints, see \cite{mystery}.

It is always possible to perform blow-ups of sufficiently small
sizes, because, by the Darboux theorem, there always exist
Darboux charts, and these contain balls.

\subsection*{Uniqueness of blow-ups of $\CP^2$}

A manifold obtained from $\tCP$ by a sequence of symplectic blow-ups 
is determined up to symplectomorphism by the sizes of the blow-ups. 
This is a non-trivial result; it is due to Biran and McDuff,
see \cite{McD:isotopy}.
We will not use it in proving our theorem.

\section{Equivariant symplectic blow-ups}

\subsection*{Equivariant symplectic blow-ups}

The standard action of the unitary group $U(n)$
descends to the standard blow-ups of $\C^n$.
Let $M$ be a $2n$-dimensional symplectic manifold with an action of a 
compact group $G$.  Let $\Omega \subset \C^n$ be an open subset
that contains a ball about the origin of radius greater than $r$.
Let $i \colon \Omega \to M$ be a symplectomorphism onto an
open subset of $M$ that is $G$-equivariant where $G$ acts on $\Omega$
through some homomorphism $G \to U(n)$.  Then the $G$-action
naturally extends to the symplectic blow-up of $M$ obtained from $i$.
If the action on $M$ is Hamiltonian, its moment map naturally extends
to the blow-up.  (This follows from the fact that the $S^1$-action
by which we collapse is the centre of $U(n)$.)
When $G \cong (S^1)^n$ we call this a \emph{toric blow-up}.

If $G$ is Abelian then, after possibly composing the inclusion map 
$i \colon \Omega \to M$ on the right by an element of $U(n)$,
we may assume that $G$ acts on $\C^n$ by rotations of the coordinates.
Specifically, if the action is toric, we may assume that $G$ acts 
on $\C^n$ through an isomorphism with $(S^1)^n$, and if $G=S^1$,
we may assume that its action on $\C^n$ is
$$ \lambda \cdot (z_1,\ldots,z_n) 
   = (\lambda^{m_1} z_1, \ldots, \lambda^{m_n} z_n)$$
where $m_1,\ldots,m_n$ are integers. These integers are called 
the \emph{isotropy weights} at $p$; they are determined, up to 
permutation, by the $S^1$-action and the symplectic form.

\subsection*{Size and rational length}

The \emph{size}, or \emph{rational length}, 
of an interval $AB$ of rational slope in $\R^n$ 
is the largest positive number $\lambda$
such that the vector $\frac{1}{\lambda}\overrightarrow{AB}$ 
is in the lattice $\Z^n$.
This is compatible with the previous notions of ``size"
(normalized symplectic area):
in a symplectic toric manifold $(M_\Delta,\omega_\Delta,\Phi_{\Delta})$
with moment map image $\Delta$, the preimage of an edge of $\Delta$ 
is a 2-sphere in $M_\Delta$ whose normalized symplectic area 
is equal to the rational length of the edge.

\subsection*{Toric blow-ups in dimension 4}

\begin{figure}
\setlength{\unitlength}{0.0006in}
\begingroup\makeatletter\ifx\SetFigFont\undefined%
\gdef\SetFigFont#1#2#3#4#5{%
  \reset@font\fontsize{#1}{#2pt}%
  \fontfamily{#3}\fontseries{#4}\fontshape{#5}%
  \selectfont}%
\fi\endgroup%
{\renewcommand{\dashlinestretch}{30}
\begin{picture}(4199,3300)(0,-10)
\dashline{60.000}(450,1050)(450,150)(1350,150)
\path(480.000,2905.000)(450.000,3025.000)(420.000,2905.000)
\path(450,3025)(450,1050)
\path(3405.000,120.000)(3525.000,150.000)(3405.000,180.000)
\path(3525,150)(1350,150)
\put(975,1350){\makebox(0,0)[lb]{\smash{{{\SetFigFont{12}{14.4}{\rmdefault}{\mddefault}{\updefault} $\sss |z_1|^2+|z_2|^2 \geq r^2$ }}}}}
\path(450,1050)(1350,150)
\put(225,525){\makebox(0,0)[lb]{\smash{{{\SetFigFont{12}{14.4}{\rmdefault}{\mddefault}{\updefault} $\sss \epsilon$ }}}}}
\put(0,3150){\makebox(0,0)[lb]{\smash{{{\SetFigFont{12}{14.4}{\rmdefault}{\mddefault}{\updefault} $\sss |z_2|^2/2$ }}}}}
\put(825,0){\makebox(0,0)[lb]{\smash{{{\SetFigFont{12}{14.4}{\rmdefault}{\mddefault}{\updefault} $\sss \epsilon$ }}}}}
\put(900,600){\makebox(0,0)[lb]{\smash{{{\SetFigFont{12}{14.4}{\rmdefault}{\mddefault}{\updefault} $\sss \epsilon$ }}}}}
\put(3600,75){\makebox(0,0)[lb]{\smash{{{\SetFigFont{12}{14.4}{\rmdefault}{\mddefault}{\updefault} $\sss |z_1|^2/2$ }}}}}
\end{picture}
}
\caption{Blow-up of $\C^2$ of size $\displaystyle \eps = \frac{r^2}{2}$}
\labell{fig:BlC2}
\end{figure}

The moment map image of the standard symplectic blow-up of $\C^2$
of size $\eps$ is shown in Figure \ref{fig:BlC2}. It is obtained
from the moment map image $\R_+^2$ of $\C^2$ by ``chopping off" a triangle
whose sides have size $\eps$.

An equivariant blow-up of size $\eps$ of a symplectic toric manifold
amounts to ``chopping off a corner of size $\eps$" of its moment map image
$\Delta$.  This can be done if and only if there exist two adjacent edges
in $\Delta$ whose sizes are both strictly greater than $\eps$.

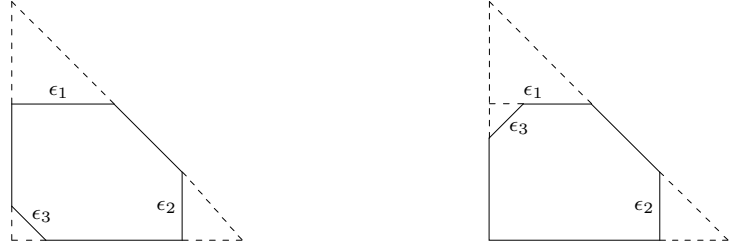
\begin{figure}
\setlength{\unitlength}{0.0006in}
\begingroup\makeatletter\ifx\SetFigFont\undefined%
\gdef\SetFigFont#1#2#3#4#5{%
  \reset@font\fontsize{#1}{#2pt}%
  \fontfamily{#3}\fontseries{#4}\fontshape{#5}%
  \selectfont}%
\fi\endgroup%
{\renewcommand{\dashlinestretch}{30}
\begin{picture}(6537,2139)(0,-10)
\path(225,1212)(1125,1212)
\path(1725,612)(1725,12)
\path(225,312)(525,12)
\path(225,1212)(225,312)
\path(1125,1212)(1725,612)
\path(1725,12)(525,12)
\path(4725,1212)(4425,912)
\path(5925,612)(5925,12)
\path(4725,1212)(5325,1212)(5925,612)
\path(5925,12)(4425,12)(4425,912)
\dashline{60.000}(225,1212)(225,2112)(1125,1212)
\dashline{60.000}(1725,612)(2250,12)(1725,12)
\dashline{60.000}(525,12)(225,12)(225,312)
\dashline{60.000}(4425,912)(4425,2112)(5325,1212)
\dashline{60.000}(5925,612)(6525,12)(5925,12)
\dashline{60.000}(4425,1212)(4725,1212)
\put(550,1287){\makebox(0,0)[lb]{\smash{{{\SetFigFont{12}{14.4}{\rmdefault}{\mddefault}{\updefault}$\sss \epsilon_1$}}}}}
\put(1500,280){\makebox(0,0)[lb]{\smash{{{\SetFigFont{12}{14.4}{\rmdefault}{\mddefault}{\updefault}$\sss \epsilon_2$}}}}}
\put(400,200){\makebox(0,0)[lb]{\smash{{{\SetFigFont{12}{14.4}{\rmdefault}{\mddefault}{\updefault}$\sss \epsilon_3$}}}}}
\put(4725,1287){\makebox(0,0)[lb]{\smash{{{\SetFigFont{12}{14.4}{\rmdefault}{\mddefault}{\updefault}$\sss \epsilon_1$}}}}}
\put(5700,280){\makebox(0,0)[lb]{\smash{{{\SetFigFont{12}{14.4}{\rmdefault}{\mddefault}{\updefault}$\sss \epsilon_2$}}}}}
\put(4600,970){\makebox(0,0)[lb]{\smash{{{\SetFigFont{12}{14.4}{\rmdefault}{\mddefault}{\updefault}$\sss \epsilon_3$}}}}}
%
%
\end{picture}
}
\caption{Equivariant blow-ups of sizes $\eps_1,\eps_2,\eps_3$ of $\CP^2$}
\label{fig:three_blowups}
\end{figure}

\begin{Example}
Figure \ref{fig:three_blowups} shows the moment map images 
of equivariant symplectic blow-ups of $\CP^2$ of sizes
$\eps_1,\eps_2,\eps_3$.   By the ``uniqueness of blow-ups" of $\CP^2$,
the resulting manifolds are (non-equivariantly) symplectomorphic.
\end{Example}

\subsection*{Equal toric blow-ups of $\mathbf{CP^2}$}

\begin{Lemma} \labell{T equivariant cor}
$(\CP^2,\omega_\FS)$ admits a toric blow-up of size $\eps > 0$ 
if and only if $\eps < 1$.
\ $(\CP^2,\omega_\FS)$ admits two or three toric blow-ups 
of size $\eps > 0$ if and only if $\eps < \half$.
\ $(\CP^2,\omega_\FS)$ does not admit four or more toric blow-ups 
of equal sizes.
\end{Lemma}

\begin{proof}
By Lemma \ref{toric is standard},
the moment map image of $\CP^2$ is a triangle in which all edges have size $1$.
So a toric blow-up of size $\eps > 0$ can be performed if and only if 
$\eps < 1$. 

After one such blow-up, the moment map image has edges of sizes 
$1-\eps$, $\eps$, $1 - \eps$, $1$ (ordered cyclically), as shown in Figure
\ref{fig:one blowup}.
If $\eps \geq \half$, no two adjacent edges have size greater than $\eps$,
so one cannot perform a second toric blow-up of size $\eps$.
If $\eps < \half$, one can perform toric blow-ups of size $\eps$
at one or both of the two endpoints of the edge of size 1.  
The resulting moment map images are shown in Figure \ref{fig:two blowups}.

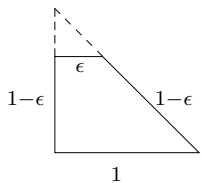
\begin{figure}
\setlength{\unitlength}{0.00083333in}
\begingroup\makeatletter\ifx\SetFigFont\undefined%
\gdef\SetFigFont#1#2#3#4#5{%
  \reset@font\fontsize{#1}{#2pt}%
  \fontfamily{#3}\fontseries{#4}\fontshape{#5}%
  \selectfont}%
\fi\endgroup%
{\renewcommand{\dashlinestretch}{30}
\begin{picture}(1212,1152)(0,-10)
\dashline{60.000}(300,825)(300,1125)(600,825)
\put(430,725){\makebox(0,0)[lb]{\smash{{{\SetFigFont{12}{14.4}{\rmdefault}{\mddefault}{\updefault}$\sss \eps$}}}}}
\path(300,825)(600,825)(1200,225)
	(300,225)(300,825)
\put(925,525){\makebox(0,0)[lb]{\smash{{{\SetFigFont{12}{14.4}{\rmdefault}{\mddefault}{\updefault}$\sss 1-\eps$}}}}}
\put(650,50){\makebox(0,0)[lb]{\smash{{{\SetFigFont{12}{14.4}{\rmdefault}{\mddefault}{\updefault}$\sss 1$}}}}}
\put(0,525){\makebox(0,0)[lb]{\smash{{{\SetFigFont{12}{14.4}{\rmdefault}{\mddefault}{\updefault}$\sss 1-\eps$}}}}}
\end{picture}
}
\caption{One toric blow-up of size $\eps$ of $\CP^2$}
\labell{fig:one blowup}
\end{figure}

\begin{figure}
\setlength{\unitlength}{0.00083333in}
\begingroup\makeatletter\ifx\SetFigFont\undefined%
\gdef\SetFigFont#1#2#3#4#5{%
  \reset@font\fontsize{#1}{#2pt}%
  \fontfamily{#3}\fontseries{#4}\fontshape{#5}%
  \selectfont}%
\fi\endgroup%
{\renewcommand{\dashlinestretch}{30}
\begin{picture}(3012,1152)(0,-10)
\path(300,825)(600,825)(900,525)
	(900,225)(300,225)(300,825)
\dashline{60.000}(300,825)(300,1125)(600,825)
\dashline{60.000}(900,525)(1200,225)(900,225)
\put(450,725){\makebox(0,0)[lb]{\smash{{{\SetFigFont{12}{14.4}{\rmdefault}{\mddefault}{\updefault}$\sss \eps$}}}}}
\put(0,450){\makebox(0,0)[lb]{\smash{{{\SetFigFont{12}{14.4}{\rmdefault}{\mddefault}{\updefault}$\sss 1-\eps$}}}}}
\put(800,350){\makebox(0,0)[lb]{\smash{{{\SetFigFont{12}{14.4}{\rmdefault}{\mddefault}{\updefault}$\sss \eps$}}}}}
\put(800,675){\makebox(0,0)[lb]{\smash{{{\SetFigFont{12}{14.4}{\rmdefault}{\mddefault}{\updefault}$\sss 1-2\eps$}}}}}
\put(500,75){\makebox(0,0)[lb]{\smash{{{\SetFigFont{12}{14.4}{\rmdefault}{\mddefault}{\updefault}$\sss 1-\eps$}}}}}
\path(2100,825)(2400,825)(2700,525)
	(2700,225)(2400,225)(2100,525)(2100,825)
\dashline{60.000}(2100,825)(2100,1125)(2400,825)
\dashline{60.000}(2700,525)(3000,225)(2700,225)
\dashline{60.000}(2400,225)(2100,225)(2100,525)
\put(2250,725){\makebox(0,0)[lb]{\smash{{{\SetFigFont{12}{14.4}{\rmdefault}{\mddefault}{\updefault}$\sss \eps$}}}}}
\put(2575,675){\makebox(0,0)[lb]{\smash{{{\SetFigFont{12}{14.4}{\rmdefault}{\mddefault}{\updefault}$\sss 1-2\eps$}}}}}
\put(2600,350){\makebox(0,0)[lb]{\smash{{{\SetFigFont{12}{14.4}{\rmdefault}{\mddefault}{\updefault}$\sss \eps$}}}}}
\put(2400,75){\makebox(0,0)[lb]{\smash{{{\SetFigFont{12}{14.4}{\rmdefault}{\mddefault}{\updefault}$\sss 1-2\eps$}}}}}
\put(1725,600){\makebox(0,0)[lb]{\smash{{{\SetFigFont{12}{14.4}{\rmdefault}{\mddefault}{\updefault}$\sss 1-2\eps$}}}}}
\put(2300,350){\makebox(0,0)[lb]{\smash{{{\SetFigFont{12}{14.4}{\rmdefault}{\mddefault}{\updefault}$\sss \eps$}}}}}
\end{picture}
}
\caption{Two or three toric blow-ups of size $\eps$ of $\CP^2$}
\labell{fig:two blowups}
\end{figure}
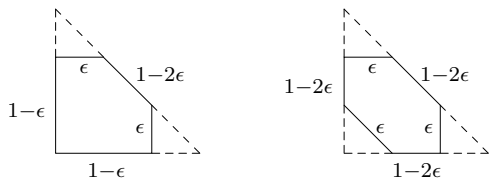

After three toric blow-ups of size $\eps$, one out of any two adjacent edges
of the moment map image has size $\eps$, see Figure \ref{fig:two blowups}.
So one cannot perform another toric blow-up of size $\eps$.
\end{proof}

\subsection*{$\mathbf{S^1}$-equivariant blow-ups in dimension 4}

\begin{Lemma} \labell{constraints}
A Hamiltonian $S^1$-manifold $(M,\omega,\Phi)$ 
admits an $S^1$-equivariant blow-up of size $\eps > 0$ centred at 
a fixed point $p \in M^{S^1}$ 
if and only if the following conditions are satisfied.
\begin{enumerate}
\item
If $p$ belongs to a $\bfZ_k$-sphere $C$ 
then $\frac{1}{2\pi} \int_C\omega > \eps$.
\item
If $p$ belongs to a two dimensional component $F$ of the fixed point set 
$M^{S^1}$ then $\frac{1}{2\pi} \int_F\omega > \eps$.
\item
$\Phi(p) + \eps < \max\limits_{m \in M} \Phi(m)$ \ and \  
$\Phi(p) - \eps > \min\limits_{m \in M} \Phi(m)$.
\item
If $p$ is an isolated fixed point and 
is a minimum or maximum for the moment map,
and if $q$ is any fixed point other than $p$, then 
$| \Phi(p) - \Phi(q) | > \eps $.
\end{enumerate}
\end{Lemma}

\begin{proof}[Proof that the conditions are necessary]
Let $\Omega \subset \C^2$ be an open set containing a ball of radius $>r$
about the origin, where $\eps = r^2/2 $, and let $i \colon \Omega \to M$ 
be a symplectic embedding that is equivariant with respect to a homomorphism
$S^1 \to U(2)$ and such that $i(0) = p$.  Let $m$ and $n$ be the isotropy
weights at $p$, so we may assume that the $S^1$-action on $\C^2$ is
$$ \lambda \cdot (z_1,z_2) = (\lambda^m z_1 , \lambda^n z_2),$$
and the moment map on $\C^2$ is
\begin{equation} \labell{Phi C2}
\Phi_{\C^2} (z) = \Phi(p) + m \frac{|z_1|^2}{2} + n \frac{|z_2|^2}{2}.
\end{equation}

The point $(z_1,0)$, for $z_1 \neq 0$, is fixed if and only if $m=0$,
and it has stabilizer $\bfZ_k$ if and only if $k = |m| \geq 2$.
A similar statement holds for the point $(0,z_2)$.
It follows that if $p$ belongs to a $\bfZ_k$-sphere $C$
or to a two dimensional component $F$ of $M^{S^1}$
then $i\inv(C)$ or $i\inv(F)$ is equal to the intersection of $\Omega$
with one of the coordinate planes. Because this intersection
contains a disk of area $\pi r^2 = 2 \pi \eps$, the size of $C$,
or $F$, must be greater than $\eps$.  This immediately implies
Conditions (1) and (2).

Condition (3) follows from the fact that the moment map for the $S^1$-action
on $\C^2$
satisfies $\sup\limits_{z \in \Omega} \Phi_{\C^2}(z) > \Phi(p) + \eps$
\, and \, $\inf\limits_{z \in \Omega} \Phi_{\C^2}(z) < \Phi(p) - \eps$.

The \emph{Duistermaat-Heckman function} of a Hamiltonian $T$-manifold
$(M,\omega,\Phi)$ is the function $\tDH_M \colon {\t}^{*} \to \R$ given by
\begin{equation} \labell{dhmd}
\tDH_M(\alpha) = 
\frac{1}{(2\pi)^k} \int_{\Phi\inv(\alpha)/T} \omega_\alpha^k / k!
\end{equation}
where $k = \half \dim M - \dim T$, and $\omega_\alpha$ is the induced symplectic form on the reduced space
$\Phi\inv(\alpha) / T$.  
The Duistermaat-Heckman function for the circle action on $\C^2$ with 
moment map \eqref{Phi C2}, when $m$ and $n$ are both positive,  is the function 
$\tDH_{\C^2} \colon \R \to \R$ given by 
$$ \tDH_{\C^2}(\alpha) = \frac{1}{mn} S(\alpha - \Phi(p)) $$
where
$$ S(t) = \begin{cases} t & \text{ if } t \geq 0 \\ 0 
                          & \text{ if } t \leq 0.  \end{cases} $$
By the Guillemin-Lerman-Sternberg formula (see \cite[Section 3.3]{GLS}
and \cite[formula (4.27)]{GGK}),
if every fixed point $q$ such that $\Phi(q) < \alpha$ is isolated,
then
\begin{equation} \labell{eq:GLS}
   \tDH_M(\alpha) = \sum_{q \in M^{S^1}} 
   \frac{1}{m_q n_q} S(\alpha - \Phi(q))  
\end{equation}
where $m_q$ and $n_q$ are the isotropy weights at $q$.

Now suppose that $\Phi$ attains its minimum at an isolated fixed point $p$.
By \eqref{Phi C2}, the set
$$ \{ z \ | \ |\Phi_{\C^2}(z) - \Phi(p)| \leq \eps \} $$
is contained in the ball of radius $r$ where $r^2/2 = \eps$,
and this ball is contained in $\Omega$, thus by (\ref{dhmd}),
\begin{equation} \labell{bigger}
 \tDH_M(\alpha) \geq \tDH_{\C^2} (\alpha) = \frac{1}{mn} S(\alpha - \Phi(p))
\end{equation}
for all $\Phi(p) \leq \alpha \leq \Phi(p) + \eps$.  
If there exists a fixed point $q$ such that $\Phi(p) < \Phi(q) < \alpha$
then, by \eqref{eq:GLS}, and since a non isolated fixed point can only occur at the maximum of $\Phi$,
\begin{equation} \labell{smaller}
 \tDH_M(\alpha) < \frac{1}{mn} S(\alpha - \Phi(p)) .
\end{equation}

Condition (4), for a minimum, follows from \eqref{bigger}
and \eqref{smaller}. Condition (4) for a maximum is proved
in the same way.
\end{proof}

\begin{proof}[Proof that the conditions are sufficient]
See \cite[Proposition 7.2]{karshon}.
\end{proof}

If $\tM'$ is obtained from $\tM$ by an equivariant blow-up, 
and $q$ is an isolated fixed point in $\tM$ different from
the point at which the blow-up is centred, then $q$ can also
be considered as a point of $\tM'$.
This follows from the fact that the image of an equivariant open 
embedding $i \colon \Omega \to \tM$, for $\Omega \subset \C^n$,
cannot contain an isolated fixed point except perhaps $i(0)$.

\begin{Lemma} \labell{cannot iterate}
Let $(M,\omega,\Phi)$ be a four dimensional compact Hamiltonian
$S^1$-manifold.  Let $(\tM,\tomega,\tPhi)$ be 
an equivariant blow-up of $(M,\omega,\Phi)$ of size $\eps >0 $ 
centred at an isolated fixed point $p \in M^{S^1}$.
Let $E \subset \tM$ be resulting exceptional divisor.
Then $\tM$ does not admit an equivariant symplectic blow-up
of size $\eps$ that is centred at a fixed point in $E$.

Moreover, suppose that $q \in E$ is an isolated fixed point.
Let $\tM'$ be obtained from $\tM$ by a sequence of equivariant 
symplectic blow-ups at points other than $q$.  Then $\tM'$ 
does not admit an equivariant symplectic blow-up of size $\eps$ 
that is centred at $q$.
\end{Lemma}

\begin{proof}
Let $m$ and $n$ be the isotropy weights at $p$.
Because $p$ is an isolated fixed point, $m$ and $n$ are non-zero.
Because the action is effective, $m$ and $n$ are relatively prime.
The exceptional divisor $E$ is symplectomorphic to $(\CP^1,\eps \omega_\FS)$
with the $S^1$-action $\lambda \cdot [z,w] = [\lambda^m z , \lambda^n w ]$.

If $k:= |m-n| \geq 2$, then $E$ is a $\bfZ_k$-sphere.
The lemma then follows from Part (1) of Lemma \ref{constraints}.

If $m=n=1$ or $m=n=-1$, then $E \subset \tM$ is fixed under the $S^1$-action.
The lemma then follows from Part (2) of Lemma \ref{constraints}.

If $|m-n| = 1$, $m$ and $n$ must have the same sign. The isolated 
fixed point $p$ is then a minimum or maximum for the moment map $\Phi$.
\ $E$ then connects the minimum or maximum of $\tPhi$ with an interior
fixed point. The lemma then follows from Parts (3) and (4) of Lemma 
\ref{constraints}.

Now suppose that $q \in E$ is an isolated fixed point
and that $\tM'$ is obtained from $\tM$ by a sequence of
equivariant symplectic blow-ups.  If $k:= |m-n| \geq 2$,
so that $E \subset \tM$ is a $\bfZ_k$-sphere of size $\eps$
passing through $q$, the additional blow-ups transform $E$
into a $\bfZ_k$-sphere through $q$ of size $\leq \eps$.
If $m=n=1$ or $m=n=-1$, 

then all the points of $E$ are fixed points so the fixed point $q$ is not isolated.
If $|m-n|=1$, then either $\tPhi(q) + \eps = \max_{m \in \tM} \tPhi(m)$
or $\tPhi(q) - \eps = \min_{m \in \tM} \tPhi(m)$.
After blowing up, the minimal value of the moment map does not
decrease and the maximal value does not increase, so in $\tM'$ we have
$\tPhi'(q) + \eps \geq \max_{m \in \tM'} \tPhi'(m)$
or $\tPhi'(q) - \eps \leq \min_{m \in \tM'} \tPhi'(m)$.
As before, in each of these cases Lemma \ref{constraints} implies 
that $\tM'$ does not admit an equivariant blow-up of size $\eps$
centred at $q$.
\end{proof}

\subsection*{Equal $\mathbf{S^1}$-equivariant blow-ups of $\mathbf{CP^2}$}

\begin{Lemma} \labell{S1 equivariant cor}
$(\CP^2,\omega_\FS)$ admits $k$ \, $S^1$-equivariant blow-ups
of size $\eps > 0$ if and only if $(k-1) \eps < 1$. 
\end{Lemma}

\begin{proof}
In Figure \ref{fig:S1 blowup} on the left
we see the graph corresponding to a $\CP^2$ of size $1$ 
with the $S^1$-action 
$$ a \cdot [z_0,z_1,z_2] = [z_0,z_1,a z_2] .$$

Suppose that $(k-1) \epsilon < 1$.
Then we can perform one blow-up at the isolated fixed point 
and $(k-1)$ blow-ups at points on the fixed surface.
(See Lemma \ref{constraints}.) 
The resulting graph is shown in Figure \ref{fig:S1 blowup}
on the right.

{
\begin{figure}
\setlength{\unitlength}{0.0006in}
\begingroup\makeatletter\ifx\SetFigFont\undefined%
\gdef\SetFigFont#1#2#3#4#5{%
  \reset@font\fontsize{#1}{#2pt}%
  \fontfamily{#3}\fontseries{#4}\fontshape{#5}%
  \selectfont}%
\fi\endgroup%
{\renewcommand{\dashlinestretch}{30}
\begin{picture}(8130,4227)(0,-10)
\put(308,3529){\blacken\ellipse{50}{50}}
\put(308,129){\blacken\ellipse{600}{100}}
\put(4508,129){\blacken\ellipse{450}{100}}
\put(4208,729){\blacken\ellipse{50}{50}}
\put(4508,729){\blacken\ellipse{50}{50}}
\put(4808,729){\blacken\ellipse{50}{50}}
\put(4508,2929){\blacken\ellipse{250}{100}}
\put(533,3454){\makebox(0,0)[lb]{\smash{{{\SetFigFont{12}{14.4}{\rmdefault}{\mddefault}{\updefault}$\sss\Phi=\alpha+1$}}}}}
\put(758,54){\makebox(0,0)[lb]{\smash{{{\SetFigFont{12}{14.4}{\rmdefault}{\mddefault}{\updefault}$\begin{subarray}{l}\Phi=\alpha\\ \text{genus}=0\\ \text{size}=1\end{subarray}$}}}}}
\put(4808,54){\makebox(0,0)[lb]{\smash{{{\SetFigFont{12}{14.4}{\rmdefault}{\mddefault}{\updefault}$\begin{subarray}{l}\Phi=\alpha\\ \text{genus}=0\\ \text{size}=1-(k-1)\epsilon\end{subarray}$}}}}}
\put(5108,654){\makebox(0,0)[lb]{\smash{{{\SetFigFont{12}{14.4}{\rmdefault}{\mddefault}{\updefault}$\sss\Phi=\alpha+\epsilon$}}}}}
\put(4808,2854){\makebox(0,0)[lb]{\smash{{{\SetFigFont{12}{14.4}{\rmdefault}{\mddefault}{\updefault}$\begin{subarray}{l}\Phi=\alpha+1-\epsilon\\ \text{genus}=0\\ \text{size}=\epsilon\end{subarray}$}}}}}
\end{picture}
}
\caption{$S^1$-equivariant blow-ups of $\CP^2$}
\labell{fig:S1 blowup}
\end{figure}
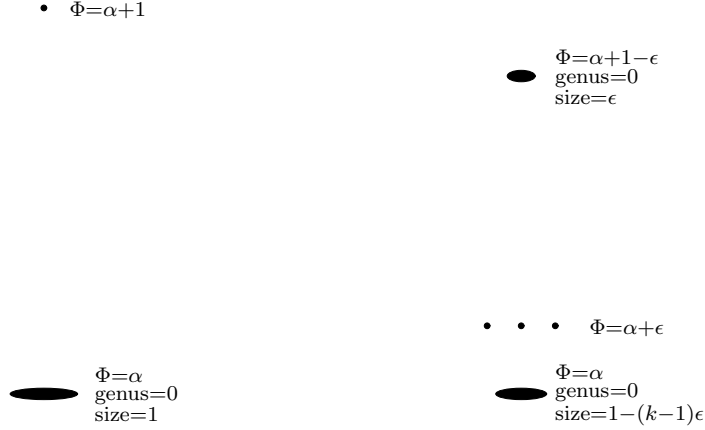
}

From the constraints described in Lemma \ref{constraints} 
it follows that the condition $(k-1) \eps < 1$ is necessary
for there to exist $k$ equivariant blow-ups of sizes $\eps > 0$.

A similar argument holds if we start from $\CP^2$ with the opposite 
$S^1$-action, $a \cdot [z_0,z_1,z_2] = [z_0,z_1,a\inv z_2] $.

Now suppose that we start with $\CP^2$ with any other $S^1$-action.
By Lemma \ref{S1 is standard}, we may assume that the action has
the form
$ \lambda \cdot [z_0,z_1,z_2] 
   = [ \lambda^a z_0 , \lambda^b z_1 , \lambda^c z_2] $
for some $a,b,c \in \Z$.  After possibly permuting the homogeneous
coordinates and multiplying them by the same power of $\lambda$,
we may assume that the action is
\begin{equation} \labell{S1 on CP2}
   \lambda \cdot [z_0,z_1,z_2] = [z_0 , \lambda^m z_1 , \lambda^n z_2] 
\end{equation}
with $0 \leq m \leq n$.  Because the action is effective, $m$ and $n$
are relatively prime.  We already covered the case that $m=0$ and $n=1$,
and, since $[z_0, \lambda z_1, \lambda z_2] = [ \lambda\inv z_0, z_1, z_2]$,
we also covered the case $m=n=1$.  So we may assume that $1 \leq m \leq n-1$.

To reduce the number of cases that we will need to consider,
we will work with \emph{extended graphs}. See \cite[page 33]{karshon}.
The extended graph of a Hamiltonian $S^1$-manifold is obtained from 
its graph by adding edges, labeled by the integer $1$, connecting 
vertices that correspond to interior fixed points to vertices that 
correspond to the minimum or maximum of the moment map, in such a way 
that the vertex corresponding to any fixed point will have exactly 
one edge going up from it and one edge going down from it.
(This graph keeps track of the ``gradient spheres" for a generic
metric. There exist other extended graphs, which correspond to
non-generic metrics, but we do not need these.)

For $(\CP^2,\omega_\FS)$ with the $S^1$-action \eqref{S1 on CP2},
with $1 \leq m \leq n-1$, and assuming that the minimal moment map
value is $\Phi=0$, we get the extended graph shown in Figure
\ref{fig:ext-graph1}, where $\ell = m-n$.

\begin{figure}
\setlength{\unitlength}{0.00083333in}
\begingroup\makeatletter\ifx\SetFigFont\undefined%
\gdef\SetFigFont#1#2#3#4#5{%
  \reset@font\fontsize{#1}{#2pt}%
  \fontfamily{#3}\fontseries{#4}\fontshape{#5}%
  \selectfont}%
\fi\endgroup%
{\renewcommand{\dashlinestretch}{30}
\begin{picture}(1192,3402)(0,-10)
\path(900,279)(300,1479)(900,3279)(900,279)
\put(900,279){\blacken\ellipse{50}{50}}
\put(300,1479){\blacken\ellipse{50}{50}}
\put(900,3279){\blacken\ellipse{50}{50}}
\put(900,150){\makebox(0,0)[lb]{\smash{{{\SetFigFont{12}{14.4}{\rmdefault}{\mddefault}{\updefault}$\sss \Phi=0$}}}}}
\put(-25,1404){\makebox(0,0)[lb]{\smash{{{\SetFigFont{12}{14.4}{\rmdefault}{\mddefault}{\updefault}$\sss \Phi=m$}}}}}
\put(950,3279){\makebox(0,0)[lb]{\smash{{{\SetFigFont{12}{14.4}{\rmdefault}{\mddefault}{\updefault}$\sss \Phi=n$}}}}}
\put(975,1779){\makebox(0,0)[lb]{\smash{{{\SetFigFont{12}{14.4}{\rmdefault}{\mddefault}{\updefault}$\sss n$}}}}}
\put(450,2379){\makebox(0,0)[lb]{\smash{{{\SetFigFont{12}{14.4}{\rmdefault}{\mddefault}{\updefault}$\sss \ell$}}}}}
\put(450,654){\makebox(0,0)[lb]{\smash{{{\SetFigFont{12}{14.4}{\rmdefault}{\mddefault}{\updefault}$\sss m$}}}}}
\end{picture}
}
\caption{Extended graph for $S^1$-action on $\CP^2$}
\labell{fig:ext-graph1}
\end{figure}
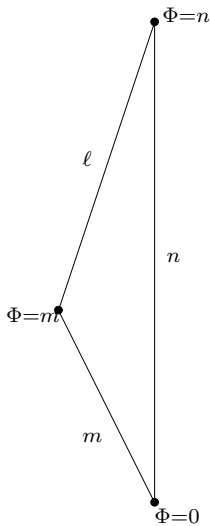

By Lemma \ref{constraints}, if we perform equivariant blow-ups
at two of the three fixed points, the sum of the sizes
of the blow-ups is less than $1$.
Fix $\eps > 0$.  By Lemma \ref{cannot iterate}, 
We can perform equivariant blow-ups of size $\eps$
only at the original three fixed points.
Hence, if $k$ is the number of blow-ups, then either
$k=1$ and $\eps < 1$, or $k \in \{2,3\}$ and $\eps < \half$.
It follows that $(k-1) \eps < 1$.

Figure \ref{fig:ext-graph2} shows the extended graph after we perform
three equivariant blow-ups of size $\eps$.

\begin{figure}
\setlength{\unitlength}{0.00083333in}
\begingroup\makeatletter\ifx\SetFigFont\undefined%
\gdef\SetFigFont#1#2#3#4#5{%
  \reset@font\fontsize{#1}{#2pt}%
  \fontfamily{#3}\fontseries{#4}\fontshape{#5}%
  \selectfont}%
\fi\endgroup%
{\renewcommand{\dashlinestretch}{30}
\begin{picture}(1363,2674)(0,-10)
\path(600,265)(300,865)
	(300,1615)(300,1615)(600,2515)
	(750,2215)(750,715)(600,265)
\put(825,1465){\makebox(0,0)[lb]{\smash{{{\SetFigFont{12}{14.4}{\rmdefault}{\mddefault}{\updefault}$\sss n$}}}}}
\put(300,2065){\makebox(0,0)[lb]{\smash{{{\SetFigFont{12}{14.4}{\rmdefault}{\mddefault}{\updefault}$\sss \ell$}}}}}
\put(300,340){\makebox(0,0)[lb]{\smash{{{\SetFigFont{12}{14.4}{\rmdefault}{\mddefault}{\updefault}$\sss m$}}}}}
\put(750,2365){\makebox(0,0)[lb]{\smash{{{\SetFigFont{12}{14.4}{\rmdefault}{\mddefault}{\updefault}$\sss n-\ell\,(=m)$}}}}}
\put(-330,1165){\makebox(0,0)[lb]{\smash{{{\SetFigFont{12}{14.4}{\rmdefault}{\mddefault}{\updefault}$\sss m+l\,(=n)$}}}}}
\put(700,400){\makebox(0,0)[lb]{\smash{{{\SetFigFont{12}{14.4}{\rmdefault}{\mddefault}{\updefault}$\sss n-m(=\ell)$}}}}}
\put(600,265){\blacken\ellipse{50}{50}}
\put(300,865){\blacken\ellipse{50}{50}}
\put(300,1615){\blacken\ellipse{50}{50}}
\put(600,2515){\blacken\ellipse{50}{50}}
\put(750,2215){\blacken\ellipse{50}{50}}
\put(750,715){\blacken\ellipse{50}{50}}
\put(600,2550){\makebox(0,0)[lb]{\smash{{{\SetFigFont{12}{14.4}{\rmdefault}{\mddefault}{\updefault}$\sss \Phi=n-\ell\eps$}}}}}
\put(825,2200){\makebox(0,0)[lb]{\smash{{{\SetFigFont{12}{14.4}{\rmdefault}{\mddefault}{\updefault}$\sss \Phi=n-n\eps$}}}}}
\put(825,680){\makebox(0,0)[lb]{\smash{{{\SetFigFont{12}{14.4}{\rmdefault}{\mddefault}{\updefault}$\sss \Phi=n\eps$}}}}}
\put(600,120){\makebox(0,0)[lb]{\smash{{{\SetFigFont{12}{14.4}{\rmdefault}{\mddefault}{\updefault}$\sss \Phi=m\eps$}}}}}
\put(-365,830){\makebox(0,0)[lb]{\smash{{{\SetFigFont{12}{14.4}{\rmdefault}{\mddefault}{\updefault}$\sss \Phi=m-m\eps$}}}}}
\put(-285,1580){\makebox(0,0)[lb]{\smash{{{\SetFigFont{12}{14.4}{\rmdefault}{\mddefault}{\updefault}$\sss \Phi=m+l\eps$}}}}}
\end{picture}
}
\caption{Three $S^1$-equivariant blow-ups of $\CP^2$}
\labell{fig:ext-graph2}
\end{figure}
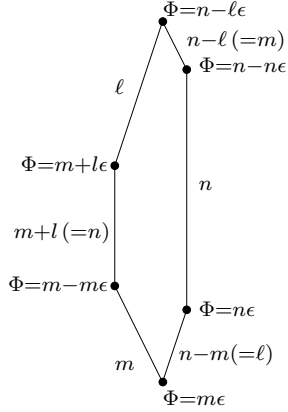
\end{proof}

\subsection*{Equivariant symplectic blow-down}

Let a compact Lie group $G$ act on a four dimensional symplectic 
manifold $(M,\omega)$.  Let $C \subset M$ be a smooth, $G$-invariant,
symplectic 2-sphere, whose homology class satisfies
$[C] \cdot [C] = -1$.
By the equivariant version of Weinstein's tubular neighborhood
theorem \cite{W}, a neighborhood of $C$ is equivariantly symplectomorphic
to a neighborhood of the exceptional divisor in a standard blow-up
of $\C^2$.
We can then equivariantly blow-down $M$ along $C$.
This yields a symplectic manifold $(N,\omega_N)$ with a $G$-action 
whose equivariant blow-up is $(M,\omega)$.

\section{Circle and torus actions on equal symplectic blow-ups of $\mathbf{CP^2}$.}

We are now ready to state our main result.  Recall that we assume that
our circle or torus actions are always effective.

\begin{Theorem} \labell{thm:main}
Let $\epsilon > 0$ be such that $1/\epsilon$ is an integer.
Let $(M_k,\omega_\epsilon)$ be a symplectic manifold
that is obtained from $(\CP^2,\omega_\FS)$ 
by $k$ simultaneous blow-ups of size $\epsilon$.

\begin{enumerate}
\item
If $k \geq 4$ then $(M_k,\omega_\epsilon)$ does not admit
a two dimensional torus action.  
\item
If $(k-1) \epsilon \geq 1$ then $(M_k,\omega_\epsilon)$
does not admit a circle action.
\end{enumerate}
\end{Theorem}

\begin{Remarks}
\begin{enumerate}
\item
The reason that Theorem \ref{thm:main} does not follow directly
from Lemmas \ref{T equivariant cor}
and \ref{S1 equivariant cor} is that a-priori
$(M_k,\omega_\epsilon)$ might admit an action which does not arise 
from an action on $\CP^2$ by performing the blow-ups equivariantly.

\item
The second author, in \cite{thesis}, showed that the conclusion 
of Theorem \ref{thm:main} holds for all $\eps \leq \frac{1}{3k 2^{2k}}$.
In joint work with Martin Pinsonnault, we further strengthened this
to all $\eps \leq 1/3$.  
Finally, Martin Pinsonnault showed that the conclusion of Theorem
\ref{thm:main} also holds for all $1/3 < \eps < 1/2$.
These results will appear elsewhere.
\end{enumerate}
\end{Remarks}

\section{Holomorphic spheres}

In this section we recall some well-known results
from Gromov's theory of holomorphic curves.  
We restrict our attention to curves of genus zero.
Wherever we do not specify the degree of smoothness, we assume
that maps are $C^\infty$ smooth and spaces of maps are equipped
with the $C^\infty$ topology.

An almost complex structure on a manifold $M$ is an automorphism
$J \colon TM \to TM$ of the tangent bundle such that 
$J^2 = -$identity.
It is \emph{tamed} by a symplectic form $\omega$
if $\omega(u,Ju) > 0$ for all $u \neq 0$.
This implies that for every embedded submanifold $C \subset M$,
if $J(TC) = TC$ then $\omega|_{TC}$ is non-degenerate.
On a symplectic manifold there always exists a tamed almost complex structure. 
If a symplectic manifold admits an action of a compact
Lie group $G$, there always exists a tamed almost 
complex structure that is $G$-invariant. 

Fix a compact symplectic manifold $(M,\omega)$.
Let $\J$ denote the space of almost complex structures $J$ on $M$
that are tamed by $\omega$.
The space $\J$ is contractible \cite{MS:intro}. 
It follows that the first Chern class of the complex vector bundle $(TM,J)$ 
is independent of the choice of $J \in \J$; we denote it $c_1(TM)$.

Fix $J \in \J$.
A \emph{parametrized $J$-sphere} is a map $f \colon \CP^1 \to M$ 
which is $J$-holomorphic, that is, which satisfies
the Cauchy-Riemann equations $df \circ j = J \circ df$.
It is called \emph{simple} if it cannot be factored through a branched 
covering of $\CP^1$. 

An embedding is a one to one immersion which is a homeomorphism with its image.
An \emph{embedded J-sphere} $C \subset M$ is the image 
of a $J$-holomorphic embedding $f \colon \CP^1 \to M$.
In particular, such a $C$ is an embedded symplectic sphere.

Let $A \in H_2(M;\Z)$ be a homology class.
For $J \in \J$, let $\calM(A,J)$ denote the set
of simple parametrized $J$-spheres in the class $A$.
Notice that if $A \cdot A = -1$ then 
every parametrized $J$-sphere in the class $A$ is simple.

Consider the universal moduli space
$$ \M(A,\J)=\{(f,J) \ \mid \ J \in \J, \ f \in \M(A,J)\}. $$
(In the notation of McDuff and Salamon, \cite[\S 3.1]{nsmall},
this space would be denoted by $\M^*(A,\Sigma;\J)$ 
with $\Sigma = \CP^1$.)

The projection 
$$ \pi \colon \M(A,\J) \rightarrow \J $$ 
is a smooth map between Fr\'echet manifolds.
A point $(f,J) \in \M(A,\J)$ is \emph{regular} for $\pi$
if and only if the differential of $\pi$ at $(f,J)$ is onto.
See \cite[Definition 3.1.4 and Remark 3.2.8]{nsmall}.

\begin{Lemma} \labell{implicit}
Let $J \in \J$.  Suppose that there exists a simple $J$-sphere $f$
such that $(f,J)$ is a regular point for $\pi$.  Then the image 
of $\pi$ contains a neighborhood $\Omega$ of $J$ in $\J$.
\end{Lemma}

\begin{proof}
Let $\J^\ell$ denote the space of almost complex structures on $M$
of type $C^\ell$ that are tamed by $\omega$, 
and let $\M(A,\J^\ell)$ denote the space of pairs $(f,J)$ 
where $J \in \J^\ell$ and where $f \colon \CP^1 \to M$
is a $J$-holomorphic map of type $C^\ell$ in the class $A$
which cannot be factored through a branched covering of $\CP^1$.
(By ``elliptic regularity", if $J \in \J^\ell$ 
and $f \colon \CP^1 \to M$ is $J$-holomorphic, then $f$ is of
type $C^\ell$. Hence, the space that we denote $\M(A,\J^\ell)$
coincides with the space that, in the notation of \cite[p.~42]{nsmall},
is denoted $\M^*(A,\Sigma;\J^\ell)$ with $\Sigma = \CP^1$.
See \cite[Proposition 3.1.9]{nsmall}.)

Let $\pi^\ell \colon \M(A,\J^\ell) \to \J^\ell$ denote 
the projection map. If $\ell$ is sufficiently large, 
$\J^\ell$ and $\M(A,J^\ell)$ are Banach manifolds
and $\pi^\ell$ is differentiable.
See \cite[Proposition 3.2.1]{nsmall}.
Because the differential of $\pi$ at $(f,J)$ is surjective
and its kernel is finite dimensional,
the differential of $\pi^\ell$ at $(f,J)$ is also surjective
and its kernel is finite dimensional.
See the last paragraph of \S 3.1 in \cite{nsmall}.
In particular, the kernel splits, that is, 
it has a closed complementary subspace in $T_{(f,J)}\M(A,\J^\ell)$.
See \cite[p.~4]{lang}.  
By the implicit function theorem for Banach spaces,
there exist a neighborhood of $(f,J)$ in $\M(A,\J^\ell)$ 
in which $\pi^\ell$ is a projection map;
in particular, the image of $\pi^\ell$ contains an open neighborhood 
$\Omega^\ell$ of $J$ in $\J^\ell$.
See \cite[Chapter I, \S 5, Cor.~2s]{lang}.
Since the $C^\ell$ topology on $\J \subset \J^\ell$ is coarser
then the $C^\infty$ topology, and $\Omega^\ell$ is open in $\J^\ell$,
the intersection $\Omega := \Omega^\ell \cap \J$ is open in $\J$.
Let $J' \in \Omega$. Since $J' \in \Omega^\ell$, there exists
a $J'$-holomorphic sphere $f \colon \CP^1 \to M$ of type $C^\ell$
in the class $A$.  Since $J'$ is smooth, by ``elliptic regularity", 
$f$ is smooth.  See \cite[Proposition 3.1.9]{nsmall}.  
The lemma follows.
\end{proof}

The group $\PSL(2,\C)$ of M\"obius transformations of $\CP^1$
naturally acts on $\calM(A,\J)$ by reparameterizations.
We say that $A$ is \emph{indecomposable} if it cannot be written 
as a sum $A = A_1 + A_2$ where $A_i \in H_2(M;\Z)$
and $\int_{A_i}{\omega} > 0$.

Recall that a continuous map between topological spaces
is called \emph{proper} if and only if the pre-image of
any compact set is compact.

\begin{Lemma} \labell{proper}
If $A$ is indecomposable then the map
\begin{equation} \labell{map on quotient}
\begin{CD}
 \M(A,\J)/\PSL(2,\C)  @>>> \J 
\end{CD}
\end{equation}
that is induced from $\pi$ is proper.
\end{Lemma}

\begin{proof}
By Gromov's compactness theorem \cite[Theorem 5.3.1]{nsmall},
since $A$ is indecomposable, if $J_n$ converges in $\J$,
then every sequence $(f_n,J_n)$ in $\M(A,\J)$ has a convergent 
sub-sequence.
Recall that a space is defined to be sequentially compact
if and only if every sequence has a convergent subsequence.
It follows that the pre-image under $\pi$ of any sequentially compact
subset of $\J$ is a sequentially compact subset of $\M(A,\J)$.
In a space whose topology has a countable bases, a subset is compact
if and only if it is sequentially compact.  
See \cite[Chap.~5, Thm.~2 and Lemma~4]{kelley}.
It follows that $\pi$ is proper.
\end{proof}

\begin{Lemma} \labell{adjunction}
Suppose that $\dim M = 4$.
Let $A \in H_2(M;\Z)$ be a homology class which is represented
by an embedded symplectic sphere $C$.  Then
\begin{enumerate}
\item
There exists an almost complex structure $J_0 \in \J$
for which $C$ is a $J_0$-sphere.
\item
For any $J \in \J$ and any parametrized $J$-sphere 
$f \colon \CP^1 \to M$ in the class $A$, the map $f$ is an embedding.
\end{enumerate}
\end{Lemma}

\begin{proof}
Construct $J_0$ as follows. Let
$$ f_0 \colon \CP^1 \to M $$
be a symplectic embedding whose image is $C$.
Define $J_0|_{TC}$ such that $f_0$ is holomorphic.
Extend it to a compatible fiberwise complex structure on the 
symplectic vector bundle $TM|_C$.  
Then extend it to a compatible almost complex structure on $M$. 
See \cite[Section 2.6]{MS:intro}.
Then $(f_0,J_0) \in \M(E,\J)$.  

By the adjunction inequality, for any $(f,J) \in \M(A,\J)$,
$$ A \cdot A - c_1(TM)(A) + 2 \geq 0,$$
with equality if and only if $f$ is an embedding.
See \cite[Cor.~E.1.7]{nsmall}.
Applying this to $(f_0,J_0)$, we get that the homology class $A$
satisfies $A \cdot A - c_1(TM)(A) + 2 = 0$.
Applying the adjunction inequality to any other $(f,J) \in \M(A,\J)$,
we get that $f$ is an embedding.
\end{proof}

\begin{Lemma} \labell{exceptional}
Let $(M,\omega)$ be a compact symplectic four-manifold.
Let $A \in H_2(M;\Z)$ be an indecomposable homology class
which is represented by an embedded symplectic sphere
and such that $c_1(TM)(A) \geq 1$.
Then for any almost complex structure $J \in \J$ there exists
an embedded $J$-sphere in the class $A$.
\end{Lemma}

\begin{proof}
Let $C \subset M$ be an embedded symplectic sphere 
such that $[C] = A$. By part 1 of Lemma \ref{adjunction}, 
$\M(A,\J)$ is non-empty. So the image of $\pi$ is a
non-empty subset of $\J$.

Let $(f,J) \in \M(E,\J)$.  By part 2 of Lemma \ref{adjunction},
$f$ is an embedding.  Since $M$ is four dimensional,
$f$ is an embedding, and $c_1(TM)(A) \geq 1$, 
by \cite{HLS}, $(f,J)$ is regular for $\pi \colon \M(A,\J) \to \J$.
By Lemma \ref{implicit} it follows that the map $\pi$ is open.
In particular, the image of $\pi$ is an open subset of $\J$.

The image of $\pi$ is equal to the image of the induced map
\eqref{map on quotient}.
Because, by Lemma \ref{proper}, this map is proper,
its image is a closed subset of $\J$.

We have shown that the image of $\pi$ is a subset of $\J$
which is nonempty, open, and closed.  Because $\J$ is connected, 
it follows that $\pi$ is onto.  This proves the Lemma.
\end{proof}

\section{Proof of the main result}

\begin{Lemma} \labell{main lemma}
Let $\epsilon > 0$ be such that $1/\epsilon$ is an integer.
Let $(M_k,\omega_\epsilon)$ be a symplectic manifold
that is obtained from $(\CP^2,\omega_\FS)$
by $k$ simultaneous blow-ups of size $\epsilon$.
Let $E_i \in H_2(M_k;\Z)$, for $i=1,\ldots,k$, be the homology classes
of the exceptional divisors. 
Let a compact connected Lie group $G$ act on $(M_k,\omega_\eps)$ 
symplectically.  Then there exist pairwise disjoint $G$-invariant 
embedded symplectic spheres $C_i \subset M_k$, for $i=1,\ldots,k$, 
whose homology classes are $[C_i] = E_i$.
\end{Lemma}

\begin{proof}
Let $C_i^0$ be the exceptional divisors in $(M_k,\omega_\eps)$,
so that $E_i = [C_i^0]$. Each $C_i^0$ is an embedded symplectic sphere.

The second homology group $H_2(M_k;\Z)$ is generated by $[\CP^1]$ 
and the $E_i$'s.
Let $N$ be the integer such that $\epsilon = \frac1N$.
Then, for each $A \in H_2(M_k;\Z)$, the size  
$\frac{1}{2\pi} \int_{A}{\omega}$ is an integer multiple
of $\frac1N$.  
It follows that any homology class whose size is equal to $\frac1N$
is indecomposable.
In particular, each of the classes $E_i$ is indecomposable.

Let $J_G$ be a $G$-invariant almost complex structure on $M_k$
tamed by $\omega_\epsilon$.
By Lemma \ref{exceptional}, for each $i$ there exists an embedded 
$J_G$-sphere $C_i \subset M_k$ such that $[C_i] = E_i$.
We now show
\begin{enumerate}
\item
$C_1,\ldots,C_k$ are disjoint;
\item
each $C_i$ is $G$-invariant;
\item
each $C_i$ is symplectic.
\end{enumerate}

For $i \neq j$, because $E_i \neq E_j$, the spheres $C_i$ and $C_j$
do not coincide.  By positivity of intersections, and since 
$E_i \cdot E_j = 0$, the spheres $C_i$ and $C_j$ are disjoint. 

Let $a \in G$.  Because $G$ is connected, $[a C_i] = [C_i] = E_i$.
By positivity of intersections and since $E_i \cdot E_i = -1$,
$aC_i$ and $C_i$ must coincide.

Because $C_i$ is an embedded $J_G$-sphere and $J_G$ is tamed by $\omega_\epsilon$,
$C_i$ is symplectic.
\end{proof}

\begin{Remark}
The assertion of Lemma \ref{main lemma} is false when $\half < \epsilon < 1$ and $k=1$. 
For a detailed example see \cite[Section 6]{thesis}.
\end{Remark}

\begin{proof}[Proof of Theorem \ref{thm:main}]
Let $\eps > 0$ be such that $1/\eps$ is an integer.
Let $(M_k,\omega_\eps)$ be a symplectic manifold
that is obtained from $(\CP^2,\omega_\FS)$
by $k$ simultaneous blow-ups of size $\eps$.

Let $G \cong (S^1)^2$  or $G \cong S^1$ act on $(M_k,\omega_\epsilon)$;
because $M_k$ is simply connected, the action is Hamiltonian.

Let $C_i$ be the spheres obtained from Lemma \ref{main lemma}.
Perform equivariant symplectic blow-downs of $(M_k,\omega_\eps)$ 
along the $C_i$'s.
By \eqref{property1} and \eqref{property2},
each such a blow-up decreases the second Betti number by one
and increases the symplectic volume by $(2\pi)^2 \epsilon^2/2$.

We get a four dimensional symplectic manifold $(N,\omega_N)$,
with an effective Hamiltonian $G$-action, 
$\dim H_2(N) = 1$ and the volume of $N$ is the same as that of 
$(\CP^2,\omega_\FS)$.
By Lemmas \ref{toric is standard} and \ref{S1 is standard},
this manifold is equivariantly symplectomorphic to $\CP^2$ 
with its standard symplectic form.

Viewing this process in reverse order, we get that $(M_k,\omega_\epsilon)$ 
is obtained from $\CP^2$ by $k$ equivariant blow-ups of size $\eps$.

The theorem now follows from Lemmas \ref{T equivariant cor}
and \ref{S1 equivariant cor}.
\end{proof}

\begin{Remark} 
Scott Baldridge \cite{baldridge1} has shown that any four-manifold which admit a symplectic structure 
and a circle action with fixed points must be rational or ruled. 
Consequently, such a manifold admits a symplectic form that is preserved by a circle action. 
Baldridge asked \cite{baldridge2} whether \emph{any} symplectic form on such a manifold is 
preserved by some circle action. Our result shows that the answer is ``no". 
For example, the manifold obtained from $\CP^2$ by four blow-ups admits circle actions and 
admits symplectic forms; by \cite{traynor}, this manifold admits a symplectic form such that each of 
the four exceptional divisors has size $1/3$; by Theorem \ref{thm:main}, this symplectic form 
is not preserved under any circle action. 
\end{Remark}

\subsection*{Acknowledgement}
Most of the work on this project has taken place when the authors
were affiliated with the Hebrew University of Jerusalem.

We are grateful to Paul Biran, Francois Lalonde, Dusa McDuff,
Martin Pinsonnault, Leonid Polterovich, and Dietmar Salamon
for valuable discussions.

\end{document}